\documentclass[12pt,a4paper]{article}
\usepackage[centertags]{amsmath}
\usepackage{amsfonts}
\usepackage{graphics}
\usepackage{graphicx}
\usepackage{amssymb}
\usepackage{amsthm}
\usepackage{pstricks,pst-node,pst-tree,pstricks-add}
\usepackage{newlfont}
\usepackage{epsfig}
\usepackage[final]{pdfpages}

\newtheoremstyle{theorem}
  {10pt}          
  {10pt}  
  {\sl}  
  {\parindent}     
  {\bf}  
  {. }    
  { }    
  {}     
\theoremstyle{theorem}
\newtheorem{theorem}{Theorem}[section]

\newtheoremstyle{defi}
  {10pt}          
  {10pt}  
  {\rm}  
  {\parindent}     
  {\bf}  
  {. }    
  { }    
  {}     
\theoremstyle{defi}

\usepackage{latexsym}
\usepackage{color}
\usepackage{amssymb}
\usepackage{amsmath}
\usepackage{amsthm} \usepackage{graphicx}
\newtheorem{lem}{Lemma}[section] \newtheorem{coro}{Corollary} 
\newtheorem{defi}{Definition}[section] \newtheorem{case}{Case}\newtheorem{rem}{Remark}[section] \newtheorem{exa}{Example}[section]
 \newtheorem{cla}{Claim} \newtheorem{ill}{Illustration}[section]

\newcommand{\bl}{\begin{lem}}
\newcommand{\el}{\end{lem}}
\newcommand{\bt}{\begin{theorem}}
\newcommand{\et}{\end{theorem}}
\newcommand{\bc}{\begin{coro}}
\newcommand{\ec}{\end{coro}}
\newcommand{\bd}{\begin{defi}}
\newcommand{\ed}{\end{defi}}
\newcommand{\bp}{\begin{proof}}
\newcommand{\ep}{\end{proof}}
\newcommand{\br}{\begin{rem}}
\newcommand{\er}{\end{rem}}
\newcommand{\be}{\begin{exa}}
\newcommand{\ee}{\end{exa}}
\newcommand{\bca}{\begin{case}}
\newcommand{\eca}{\end{case}}
\newcommand{\bcl}{\begin{cla}}
\newcommand{\ecl}{\end{cla}}
\newcommand{\bil}{\begin{ill}}
\newcommand{\eil}{\end{ill}}

\parindent 18pt\date{}
\begin{document}
\title{Perturbation of operators and approximation of spectrum}
\author{Kiran, M.N.N.Namboodiri
\footnote{Department of Mathematics,
``CUSAT'' - Cochin (INDIA) (E-mail: {\tt mnnadri@gmail.com,})},
S. Serra-Capizzano \footnote{Department 'Fisica e Matematica',Via Valleggio 11, 22100 Como (ITALY) (E-mail: {\tt stefano.serrac@uninsubria.it})}}
\maketitle
\begin{abstract}
Let $A(x)$ be a holomorphic family of bounded self-adjoint operators on a separable Hilbert space $\mathbb{H}$ and let $A(x)_n$ be the orthogonal compressions of $A(x)$ to the span of first n elements of an orthonormal basis of $\mathbb{H}$. 
The problem considered here is to approximate the spectrum of $A(x)$ using the sequence of eigenvalues of $A(x)_n$. We show that the bounds of the essential spectrum and the discrete spectral values outside the bounds of essential spectrum of $A(x)$ can be approximated uniformly on all compact subsets by the sequence of eigenvalue functions of $A(x)_n$. The known results for a bounded self-adjoint operator, are translated into the case of a holomorphic family of operators. Also an attempt is made to predict the existence of spectral gaps that may occur between the bounds of essential spectrum of $A(0)=A$ and study the effect of holomorphic perturbation of operators in the prediction of spectral gaps. As an example, gap issues of some block Toeplitz-Laurent operators are discussed. The pure linear algebraic approach is the main
 advantage of the results here.
\end{abstract}
\section{Introduction}
Perturbation theory of operators incorporates a good deal of spectral theory. There are many instances in quantum mechanics, where the perturbation of operators arises. For example, the Schrodinger operator 
\begin{equation}\label{operator}
\tilde{A}(u)= -\ddot{u}+ V\cdot u
\end{equation} 
 defined on a suitable subspace of $ L^2(\mathbb{R})$: can be viewed as a perturbation of differential operator.
 Here we discuss the linear algebraic techniques used in \cite{Bot} and \cite{KSN}, under a holomorphic perturbation of the operator.

          Let $\mathbb{H}$ be a separable Hilbert space and A be a bounded self-adjoint operator defined on $\mathbb{H}$. The spectrum of $A$ will be a compact subset of the interval $[m,M]$ denoted by $\sigma(A)$. Let $\{e_1,e_2,\ldots\} $ be an orthonormal basis for $\mathbb{H}$. Consider the finite dimensional truncations of $A$, that is, $A_n = P_nAP_n$ where $P_n$ is the projection of  $\mathbb{H}$ onto the span of first n elements $\ \{e_1,e_2,\ldots ,e_n \} $ of the basis.

 Various mathematicians have done extensive research for computing  
$\sigma(A),$ $\sigma_e(A)$ and their bounds using $\sigma(A_n)$ \cite{Arv,Bot,Dav1,Dav2,Hans}. However prediction of spectral gaps and related problems using truncation method is yet to be investigated in detail, though a brief attempt in this direction has been done in \cite{ias}. In this paper, the approximation results in \cite{Bot}, are translated in to the case of a holomorphic family of operators $A(x)$. We prove that the bounds of essential spectrum and the discrete spectral values outside the bounds of essential spectrum can be approximated uniformly on all compact subsets by sequence of eigenvalue functions.

Also, some spectral gap prediction results are proved using the finite dimensional truncations. We should mention that gap related problems were studied using analytical and variational techniques, especially for Schrodinger operators with 
     different kinds of potentials. This refers to classical Borg-type theorems which characterized the periodic potentials depending on the nature of  
     spectral gaps (see \cite{Bru,Fla,Bar} and references therein and refer to  \cite{Clark,KSN} for new perspectives). Here we try for such results in the case
     of some perturbed discrete Schrodinger operators treating them as block Toeplitz-Laurent operators.

       The following is a brief account of some developments in the linear algebraic techniques to the spectral approximation problem, which will play a key role throughout this paper.
\subsection{Linear algebraic approach}
         Let $\nu,\mu $ be the lower and upper bounds of $\sigma_e(A)$ respectively, with $A$ being self-adjoint. Let 
         $\lambda^+_R(A)\leq\ldots\leq\lambda^+_2(A)\leq\lambda^+_1(A)$ be the discrete eigenvalues of $A$ lying above $\mu$ and $\ 
         \lambda^-_1(A)\leq\lambda^-_2(A)\leq \ldots\leq\lambda^-_S(A)$ be the eigenvalues of $A$ lying below $\ \nu$. Here $R$ and $S$ can be 
         infinity.
Denote by $\ \lambda_1(A_n)\geq\lambda_2(A_n)\geq\ldots\geq\lambda_n(A_n)$ the eigenvalues of $\ A_n$.
 The following result from \rm{\cite{Bot}} is of interest in our context.
  \bt \label{Bot}
For every fixed integer $k$ we have
  \begin{equation} \label{Bot1} \nonumber
     \lim_{n\rightarrow\infty} \lambda_k(A_n) =\left\{ {\begin{array}{*{20}c}
   {\lambda ^ +  _k \left( A \right),  \, \textrm{ if} \, R\, = \,\infty \,  \textrm{ or} \,1 \le k \le R,}  \\
   {\mu ,\,\,\,\,\,\,\,\,\,\,\,\,\,  \textrm{ if} \,R<\,\infty \,  \textrm{and} \,k \ge R + 1,}  \\
\end{array}} \right.
   \end{equation}
   \begin{equation} \label{Bot2} \nonumber
    \lim_{n\rightarrow\infty} \lambda_{n+1-k}(A_n) =\left\{ {\begin{array}{*{20}c}
   {\lambda ^ -  _k \left( A \right), \,  \textrm{if} \,S\, = \,\infty \,  \textrm{ or}\,1 \le k \le S,}  \\
   {\nu ,\,\,\,\,\,\,\,\,\,\,\,\,\,  \textrm{ if} \,S<\,\infty \,    \textrm{ and} \,k \ge S + 1.}  \\
\end{array}} \right.
\end{equation}
  In particular, $\lim_{k\rightarrow\infty} \lim_{n\rightarrow\infty} \lambda_{k}(A_n) =\mu$
and $\lim_{k\rightarrow\infty} \lim_{n\rightarrow\infty} \lambda_{n+1-k}(A_n) =\nu.$ \et
\br
The above results are also true if we replace $\ A_n $ by some other sequence $\ A_{1n}$ with the property that
$\ \left\| {A_{n\,}  - A_{1n} } \right\| \to 0$ as $n \to \infty $, with $\|\cdot\|$ being the spectral norm. In order to justify this, we 
need only to recall an important  inequality concerning the eigenvalues of self-adjoint matrices $A, B$
(refer e.g. to \cite{Bha1})
\begin{equation} \label{eigen-norm}
  \left| {\lambda _k \left( A \right) - \lambda _k \left( B \right)} \right| \le \left\| {A - B} \right\|.
  \end{equation}
\er  
  Now we give some definitions and two related theorems from {\rm \cite{Arv}}.
 \bd {\rm \cite{Arv}} {Essential points :}
           A real number $\ \lambda$ is an essential point if for every open set $U$ containing $\ \lambda,\, \lim_{n\rightarrow\infty}N_n(U) = \infty$, where
             $\ N_n(U)$ is the number of eigenvalues of $\ A_n$ in $U$.
 \ed
\bd {\rm \cite{Arv}} {Transient points :}
             A real number $\ \lambda$ is transient if there is an open set $U$ containing $\lambda,$ such that $\sup_{n\geq1}N_n(U) < \infty$.\ed

\bd  {\rm \cite{Arv}} The degree of an operator $A$ is defined by the relation
\begin{equation} \nonumber
\textrm{ deg}(A)= \sup_{n{\geq 1}} \textrm{ rank} (P_nA - AP_n).
\end{equation} \ed
 \bd {\rm \cite{Arv}} $A$ is an operator in the Arveson's class if
 \begin{equation} \nonumber
A = \Sigma_n A_n,  \textrm{ deg}( A_n)<\infty \, \textrm{for every n and }  \,\, \Sigma_n (1+\textrm{ deg}(A_n)^{\frac{1}{2}})\|A_n\|<\infty
\end{equation}  \ed
\bt\label{Arv} \cite{Arv} If A is a bounded self-adjoint operator, and 
if we denote  
\begin{equation} \nonumber
\Lambda = \{ \lambda \in R ; \lambda = \lim \lambda_n,\lambda_n\in\sigma(A_n) 
\end{equation} 
and $\ \Lambda_e,$ the set of all essential points, then 
\begin{equation} \nonumber
\sigma(A)\subseteq\Lambda \subseteq [m,M]\,\textrm{and}\, \sigma_e(A)\subseteq\Lambda_e .
\end{equation}
 \et
\bt \cite{Arv} If $A$ is a bounded self-adjoint operator in the Arveson's class, then  $\ \sigma_e(A)= \Lambda_e $ and every point in $\ \Lambda $ is either 
transient or essential.\label{Arv1}\et
  The subsequent theorem taken from  {\rm \cite{Bot}} denies the existence of spurious eigenvalues (points in $\ \Lambda$ which are not spectral 
  values), under the assumption that the essential spectrum is connected.

\bt  \label{nogap}
If $A$ is a self-adjoint operator and if $\ \sigma_e(A)$ is connected, then $\ \sigma(A) = \Lambda $.
\et
\br It is worthwhile to notice that the connectedness of essential spectrum enables us to compute the spectrum using finite dimensional truncations. 
\er
 The paper is organized as follows. In Section 2, we describe the generalized versions of the approximation results in the case of a one parameter holomorphic family of 
 operators. We follow the definitions by Kato in \cite{Kat}. In Section 3, we prove results 
 regarding the prediction of spectral gap with some examples. Also we make observations of what happens to the spectral gaps under holomorphic perturbation. In the fourth section, we report results on the spectral gaps of some block Toeplitz-Laurent operators. We present the modified version of discrete Borg's theorem with the techniques used in \cite{KSN}. A concluding section ends the paper.
\section{Spectrum under perturbation}

Let $ A\left( x \right)$ be a holomorphic family of operators with domain $D_0$ in the 
complex plane. That is for each $y \in D_0$, the following limit 
\begin{equation}\nonumber
\lim_{x\rightarrow y}\frac{\left\|A(x)-A(y)\right\|}{\left|x-y\right|}
\end{equation}
should exist and must be finite.
Our aim is to study the changes in the behavior of spectrum, under these perturbations. 
We generalize the approximation techniques used in the case of a 
single operator, to a holomorphic family of operators. 
First we define the approximation number functions as follows.
\bd
Consider the singular number $s_k$, $k$ natural number,
\begin{equation}\nonumber
s_k(A\left( x \right))=\inf \{\left\|A\left( x \right)-F\right\|, \textrm{ rank} (F)\leq k-1\} x \in D_0
\end{equation}
is the $k^{th}$ \textbf{approximation number function} of $A\left( x \right)$.
\ed
Clearly we have for each $x \in D_0$,
\begin{equation} \label{Essential norm}
 \left\| {A\left( x \right)} \right\| = s_1 \left( {A\left( x \right)} \right) \ge s_2 \left( {A\left( x \right)} \right)\ge  \ldots \ge s_k \left( 
 {A\left( x \right)} \right) \ge \ldots  \ge 0.
\end{equation}
Recall the definition of essential norm.
\bd  $\left\|A\left( x \right)\right\|_{ess}$ = $\inf \{\left\|A\left( x \right)-K\right\| ,K \textrm{ compact}\}$, $x \in D_0$.
\ed
The following lemmas are easy consequence of the continuity of $A(x).$
\bl \label{GolP}
$ s_k(A\left(. \right)) \rightarrow \left\|A\left(. \right)\right\|_{ess} $ as $k\rightarrow \infty$ uniformly on all compact subsets of $D_0$.
\el
\bp
Consider the sequence of functions $f_k(x)=s_k(A\left( x \right)).$ From \cite{Gol}, we have for each x, 
\begin{equation} \nonumber
f_k(x)=s_k(A\left( x \right)) \rightarrow \left\|A\left( x \right)\right\|_{ess}. 
\end{equation}  
Also since
\begin{equation}
\left| {f_k(x) - f_k(y)} \right| = \left| {s_k \left( {A\left( x \right)} \right) - s_k \left( {A\left( y  \right) } \right)} \right|\le \left\| {A\left( x \right) - A\left( y \right)} \right\|,\\
\end{equation}  
and $A(x)$ is holomorphic, we observe that each functions in the sequence are continuous. Hence using the monotonicity of the sequence of functions in (\ref{Essential norm}), we conclude that the 
convergence is uniform in each compact subsets, by Dini's theoerem (see pp.150 \cite{Rud}). Hence the proof.
\ep
  Now we consider the truncations $ A\left( x \right)_n = P_nA(x)P_n$ and singular numbers
 $ s_k(A\left( x \right)_n)=$ $\inf \{\left\|A\left( x \right)_n-F_n\right\| , \textrm{ rank} (F_n)\leq k-1\}$.
 \bl\label{BotP1}
  $ s_k(A\left( x \right)_n) \rightarrow s_k(A\left( x \right)) $ as $n\rightarrow \infty$, for each $k$, 
 and the convergence is uniform on all compact subsets of $D_0$.
  \el 
 \bp  Our first observation is that the sequence of functions $ f_{n,k}(x)=s_k(A\left( x \right)_n)$ 
 form an equicontinuous family of functions. This follows from the following inequality.
 \begin{equation} \nonumber
\begin{array}{l}
 \left| {f_{n,k} \left( x \right) - f_{n,k} \left( y \right)} \right| = \left| {s_k \left( {A\left( x \right)_n } \right) - s_k \left( {A\left( y 
 \right)_n } \right)} \right| \le \left\| {A\left( x \right)_n  - A\left( y \right)_n } \right\| \\
  \le \left\| {A\left( x \right) - A\left( y \right)} \right\|.\end{array}
  \end{equation}
 Also from the interlacing theorem for singular values (see \cite{Bha1} for the proof), we have
\begin{equation}\nonumber
f_{n,k}(x)=s_k(A\left( x \right)_n) \geq s_k(A\left( x \right)_{n-1})=f_{n-1,k}(x),
\end{equation}
for each k and for every $x \in D_0.$
  Hence the sequence of singular value functions form a monotone sequence of functions.
  Also by Theorem (1.1) of \cite{Bot}, 
\begin{equation}\nonumber
f_{n,k}(x)=s_k(A\left( x \right)_n) \rightarrow s_k(A\left( x \right))\,\,
\textrm{as}\,\,n\rightarrow \infty,
  \end{equation}
for each k and for all x $\in D_0.$
 Now by Dini's theorem, the convergence is uniform on all compact subsets
 of $D_0$ and the proof is completed.
\ep
  For the rest of this paper, we assume that $A(x)$ is self-adjoint for each x.
Let $\nu(x),\mu(x) $ be the lower and upper bounds of $\sigma_e(A(x))$ respectively, and also let the numbers
$~\lambda^+_R(A(x))\leq\ldots\leq\lambda^+_2(A(x))\leq\lambda^+_1(A(x))$ be the discrete eigenvalues of $A(x)$ lying above $\mu(x)$, and $\lambda^-_1(A(x))\leq\lambda^-_2(A(x))\leq \ldots\leq\lambda^-_S(A(x))$ be the eigenvalues lying below 
$\nu(x)$. Here $R$ and $S$ can be infinity.
 The quantities $ \lambda_{1,n}(x)\geq\lambda_{2,n}(x)\geq\ldots\geq\lambda_{n,n}(x)$ denote the eigenvalues of $\ A(x)_n$ in non increasing order.
  \bt \label{perturbspectrum}
  \begin{equation}\nonumber
  \lim_{n\rightarrow\infty} \lambda_{k,n}(x) =\left\{ {\begin{array}{*{20}c}
   {\lambda ^ +  _k \left( x \right),  \, \,\textrm{ if}  \, \,R\, = \,\infty \, \, \textrm{or}  \,1 \le k \le R,}  \\
   {\mu(x) ,\,\,\,\,\,\,\,\,\,\,\,\,\,  \textrm{ if}   \,R<\,\infty \,  \textrm{and}  \,k \ge R + 1,}  \\
\end{array}} \right.
\end{equation}
\begin{equation}\nonumber
  \lim_{n\rightarrow\infty} \lambda_{n+1-k,n}(x) =\left\{ {\begin{array}{*{20}c}
   {\lambda ^ -  _k \left( x \right),  \,\,\, \textrm{if}\,\,  \,S\, = \,\infty \, \textrm{ or}  \,1 \le k \le S,}  \\
   {\nu(x) ,\,\,\,\,\,\,\,\,\,\,\,\,\,  \textrm{if} \,\,  \,S<\,\infty \, \,\textrm{and} \, \,k \ge S + 1.}  \\
\end{array}} \right.
\end{equation} 
  In particular, 
  \begin{equation}\nonumber
\lim_{k\rightarrow\infty} \lim_{n\rightarrow\infty} \lambda_{k,n}(x) =\mu(x)
\,\,\textrm{and}\,\, \lim_{k\rightarrow\infty} \lim_{n\rightarrow\infty} \lambda_{n+1-k,n}(x) =\nu(x).
\end{equation}
Furthermore, in each of the cases given above, the convergence is uniform on all compact subsets of $D_0$.
\et
\bp For each fixed $x \in D_0$, these limits exist by Theorem (\ref{Bot}). We observe the fact that 
the sequence of eigenvalue functions, $f_{n,k}(x)=\lambda _{k,n}(x)$ form an equicontinuous family of functions, from the following inequalities.
\begin{eqnarray*}
 \left| {f_{n,k}(x) - f_{n,k}(y)} \right| &=& \left| {\lambda _{k,n} \left( x \right) - \lambda _{k,n} \left( y \right)} \right|\le \left\| {A\left( x \right)_n  - A\left( y \right)_n } 
 \right\| \\
  & \le & \left\| {A\left( x \right) - A\left( y \right)} \right\|.
\end{eqnarray*}
Also by Cauchy's interlacing theorem for eigenvalues,
\begin{equation}\nonumber
\lambda_{1,n+1}(x)\geq\lambda_{1,n}(x)\geq\lambda_{2,n+1}(x)\geq\ldots\lambda_{n,n+1}(x)\geq\lambda_{n,n}(x)\geq \lambda_{n+1,n+1}(x),
\end{equation} 
for each $x \in D_0.$ In particular, for each k and for every $x \in D_0,$
\begin{equation}\nonumber
f_{n+1,k}(x)=\lambda_{k,n+1}(x)\geq\lambda_{k,n}(x)=f_{n,k}(x).
\end{equation} 
Hence $f_{n,k}(.)$ forms a monotone sequence of continuous functions that converges point wise.
Therefore by  Dini's theorem,
the convergence is uniform on all compact subsets of $D_0$.
Hence the proof is completed.
 \ep
\br
Using Theorem (\ref{perturbspectrum}), we can approximate the discrete spectrum 
of a holomorphic family of operators, lying outside the bounds of essential spectrum
 by the eigenvalue functions of truncations uniformly on all compact subsets. 
\er 
It was observed in \cite{Bot} that norm of ${A_n}^{-1}$ is uniformly bounded if A is invertible and the essential spectrum is connected. The perturbed version of this result is proved below.
\bc
Let $A(x)$ be a holomorphic family of bounded self-adjoint operators such that $\sigma_e(A(x))$ is connected for all x in the 
domain $D_0.$ Then 
\[
\lim_{n \rightarrow \infty} \left\|(A(x)_n-\lambda I_n)^{-1}\right\|=\left\|(A(x)-\lambda I)^{-1}\right\|\,\,\textrm{for every}\,\,\,\lambda \in \mathbb{C}-\mathbb{R}.
\]
Also the convergence is uniform on all compact subsets of $D_0.$
\ec
\bp
By Theorem (\ref{nogap}), $\sigma(A(x))= \Lambda(A(x))$. Hence we can easily observe the following.
\[
d(z,\sigma(A(x)_n))\rightarrow d(z,\Lambda(A(x))=d(z,\sigma(A(x)))\,\,\textrm{for every complex number z.}
\]
Therefore, for every non real z,
\[
\left\|(A(x)_n-\lambda I_n)^{-1}\right\|=\frac{1}{d(z,\sigma(A(x)_n))}\rightarrow \frac{1}{d(z,\sigma(A(x)))}=\left\|(A(x)-\lambda I)^{-1}\right\|.
\]
Also the convergence is uniform on all compact subsets of $D_0$ as observed in the previous theorems.
\ep
\section{Gaps in the essential spectrum}
       Now we consider the problem to locate the gaps in the essential spectrum if any. The following theorem is an attempt to predict the existence of spectral gaps, using the finite dimensional truncations.
         We use the notation $\#S$ to denote the number of elements in the set S.
\bt\label{Kb1}
 Let A be a bounded self-adjoint operator, and 
 \\$ \lambda_{n1}(A_n)\geq\lambda_{n2}(A_n)\geq...\geq\lambda_{nn}(A_n)$ be the eigenvalues 
of $A_n$ arranged in decreasing order. For each  positive integer n, let $\ \{w_{nk} : k = 1,2,... n\}$  be a set of numbers such that $\ 0 \leq 
w_{nk}\leq 1$ and
 $ \sum\limits_{k = 1}^n w_{nk}=1$. If there exists a $\delta>0$ and $ K>0$ such that
\begin{equation}\label{Average}
\# \left\{ {\lambda _{nj} ;\left| {\sum\limits_{k = 1}^n {w_{nk} \lambda _{_{nk} }  - \lambda _{nj} } } \right| < \delta } \right\} < K
\end{equation}
and in addition if $\ \sigma_e(A) $ and $\ \sigma(A) $ has the same upper and lower bounds, then $\ \sigma_e(A)$
 has a gap.
 \et
\bp
Consider the set $S=\left\{{\sum\limits_{k = 1}^n {w_{nk} } \lambda _{nk}, n=1,2,3 \ldots} \right\}$ and
 observe that $\lambda _{nn}  \le \sum\limits_{k = 1}^n {w_{nk} } \lambda _{nk}\le\lambda _{n1}$. 
Also since each $\lambda _{nj}$'s lie in the interval $[m,M],$ we see that the set S is 
contained in the interval $[m,M]=[\nu,\mu].$
\case
Assume that S is a finite set, say $S=\{a_1,a_2,a_3 \ldots a_m\}$. In this case,
the value of the sum  $\sum\limits_{k = 1}^n {w_{nk} } \lambda _{nk}$ equals some of the numbers $a_i$'s
for infinitely many n. Let $a_1,a_2,a_3 \ldots a_p$ be those numbers. 
That is 
\begin{equation}\nonumber
\sum\limits_{k = 1}^n {w_{nk} } \lambda _{nk}=a_i\,\,\,\, \textrm{for infinitely many n where }i=1,2,\ldots p.
\end{equation}
From this and by the condition (\ref{Average}), for each $i=1,2,\ldots p,$ we have 
\begin{equation}\nonumber
N_n(a_i-\delta,a_i+\delta)=\# \left\{ \lambda _{nj} ;\left| {a_i  - \lambda _{nj} } \right| < \delta  \right\} < K\,\,\,
\textrm{for infinitely many n.}
\end{equation}
Hence $N_n(a_i-\delta,a_i+\delta)$ will not go to infinity as n goes to infinity.
Therefore no number in the interval $(a_i-\delta,a_i+\delta)$ is an essential point. 
Since the essential spectrum is contained in the set of all essential points, by Theorem (\ref{Arv}),
there is no essential spectral values in this interval. 
Also since each $a_i$ lies between the bounds of essential spectrum, we can choose an appropriate 
$\epsilon>0$ such that $(a_i-\epsilon,a_i+\epsilon)$ lies between the bounds and contained in the
interval $(a_i-\delta,a_i+\delta)$. Then the interval $(a_i-\epsilon,a_i+\epsilon)$ is a spectral gap.
\case
Now we consider the case when S is an infinite set.
Hence S has a limit point in $\mathbb{R}$. Now if $\ w_0 $ is a limit point of the set S,
 then we have $\ \nu \leq w_0 \leq \mu$.
  
Now the interval
  $\ \left( {w_0  - {\raise0.7ex\hbox{$\delta $} \!\mathord{\left/
 {\vphantom {\delta  2}}\right.\kern-\nulldelimiterspace}
\!\lower0.7ex\hbox{$2$}},w_0  + {\raise0.7ex\hbox{$\delta $} \!\mathord{\left/
 {\vphantom {\delta  2}}\right.\kern-\nulldelimiterspace}
\!\lower0.7ex\hbox{$2$}}} \right)$ will contain infinitely many points from the set S.
 Corresponding to these points, there are infinitely many $A_n$'s for which the number of eigenvalues in
 $\ \left( {w_0  - {\raise0.7ex\hbox{$\delta $} \!\mathord{\left/
 {\vphantom {\delta  2}}\right.\kern-\nulldelimiterspace}
\!\lower0.7ex\hbox{$2$}},w_0  + {\raise0.7ex\hbox{$\delta $} \!\mathord{\left/
 {\vphantom {\delta  2}}\right.\kern-\nulldelimiterspace}
\!\lower0.7ex\hbox{$2$}}} \right)$ is bounded by K due to (\ref{Average}). Hence
the sequence $\ N_n \left( {w_0  - \frac{\delta }{2},w_0  + \frac{\delta }{2}} \right)$ 
will not go to infinity, since a subsequence of it, is bounded by K.
 Hence no point in the interval $\ \left( {w_0  - {\raise0.7ex\hbox{$\delta $} \!\mathord{\left/
 {\vphantom {\delta  2}}\right.\kern-\nulldelimiterspace}
\!\lower0.7ex\hbox{$2$}},w_0  + {\raise0.7ex\hbox{$\delta $} \!\mathord{\left/
 {\vphantom {\delta  2}}\right.\kern-\nulldelimiterspace}
\!\lower0.7ex\hbox{$2$}}} \right)$ is an essential point. 
Since the essential spectrum is contained the set of all essential points, by Theorem (\ref{Arv}),
$\ \left( {w_0  - {\raise0.7ex\hbox{$\delta $} \!\mathord{\left/
 {\vphantom {\delta  2}}\right.\kern-\nulldelimiterspace}
\!\lower0.7ex\hbox{$2$}},w_0  + {\raise0.7ex\hbox{$\delta $} \!\mathord{\left/
 {\vphantom {\delta  2}}\right.\kern-\nulldelimiterspace}
\!\lower0.7ex\hbox{$2$}}} \right)$ contains no essential spectral values. 
Hence, as in the case 1, we can choose an $\epsilon>0,$ such that the interval 
$(w_0-\epsilon,w_0+\epsilon)$ is a spectral gap between the bounds of the 
essential spectrum and the proof is completed.
\ep
\br There is possibility for the presence of discrete eigenvalues inside the gaps in the above case.\er
\br The special case which is more interesting is when $\ w_{nk} = \frac{1}{n},$ for all n.
In that case, we are actually looking at the averages of eigenvalues of truncations and these averages can be
computed using the trace at each level.\er
\textbf{Special Choice I}\\
Let us consider an instance where these weights $\ w_{nk}$ arises naturally associated to a self-adjoint operator on a Hilbert space.
Let $\ A_n  = \sum\limits_{k = 1}^n {\lambda _{n,k} Q_{n,k}} $ be the spectral resolution of $\ A_n$.
 Define $\ w_{nk}=\left\langle {Q_{n,k} e_1,e_1} \right\rangle.$ Then $\ 0 \leq w_{nk}\leq 1$ and $\ \sum\limits_{{\ k = 1}}^{\  n} {w_{nk}  = 1}$. Now
\begin{equation} \nonumber
 {\sum\limits_{k = 1}^n {w_{nk} } \lambda _{nk} }={\sum\limits_{k = 1}^n \lambda _{nk}{\left\langle {Q_{n,k} e_1,e_1}\right\rangle}} =\left\langle {A_n e_1,e_1} \right\rangle=\left\langle {A e_1,e_1} \right\rangle = a_{1 1}.
\end{equation}
Therefore by Theorem (\ref{Kb1}), if there exists a $\delta>0$ and a $K>0,$ such that 
\begin{equation} \nonumber
\# \left\{ {\lambda _{nj} ;\left|a_{1 1} - \lambda_{nj} \right| < \delta }\right\} < K
\end{equation}
then there exists a gap in the essential spectrum of A. Hence if the first entry in the matrix representation of A, is not an essential point, then there exists a gap in the essential spectrum.
\br
All points of the form $\left\langle {A e_j,e_i} \right\rangle = a_{i j}$ are in the numerical range which lies between the bounds of essential spectrum in the case that the bounds coincide with the bounds of the spectrum. Hence in that case, 
if $a_{i j}$ is not an essential point for some $i,j$, then that will lead to the existence of a spectral gap. That means if any one of the entries in the matrix representation of A is not an essential point, then there exists a gap in the essential spectrum as indicated in the above special choice of $\ w_{nk}.$
\er
The following is an example where the first entry $a_{1 1}$ is a transient point and the spectral gap prediction is valid.
\be
Define a bounded self-adjoint operator A on $l^2(\mathbb{N})$, as follows.
\[
A(x_n)=(x_{n-1}+x_{n+1})+(v_nx_n), x_0=0;
\]
where the periodic sequence $v_n=(1,2,3,1,2,3,\ldots).$
This is a discretized version of the well known Schrodinger operator. The matrix representation of A results the 
block Toeplitz operator with corresponding matrix valued symbol given by
\[
 \tilde{f}\left( \theta  \right) = \left[ {\begin{array}{*{20}c}
   {1} & {1}&{e^{i\theta} }  \\
   {1} & {2}&{1} \\
   {e^{-i\theta}} & {1}&{3} \\ 
\end{array}} \right].
\]
As indicated in the special choice above, by Theorem (\ref{Kb1}), 
if $\left\langle A(e_1),e_1\right\rangle=1$ is a transient point, then $\sigma_e(A)$ has a gap.
This is evident in this example, since from \cite{bottcher},
\begin{equation} \nonumber
\sigma_{\hbox{ess}}(A)= \bigcup_{j=1}^3 \left[\inf_\theta(\lambda_j(\tilde{f}\left( \theta  \right)),\sup_\theta(\lambda_j(\tilde{f}\left( \theta  \right))\right]
\end{equation} 
where $\lambda_j(\tilde{f}\left( \theta  \right))$ are the eigenvalues of $\tilde{f}\left( \theta  \right).$
A straightforward numerical computation of the eigenvalue functions gives
\begin{equation} \nonumber
\sigma_{\hbox{ess}}(A)=\left[-0.2143,0.3249\right]\cup\left[1.4608, 2.5392\right]\cup \left[3.6751, 4.2143\right].
\end{equation} 
Also since A is in the Arveson's class (all band limited matrices comes in this class), the point 1 lies in the gap, is a transient point. Hence the prediction of the existence of gap, in Theorem (\ref{Kb1}), is valid in this example.
\ee

\textbf{Special Choice II}\\
By invoking Theorem (\ref{Bot}), there exists a sequences of eigenvalues of truncations $\ \lambda _{n_l},\lambda _{n_m}$ such that
$\lim_{n_l\rightarrow\infty} \lambda_{n_l} =\nu$ and $\lim_{n_m\rightarrow\infty} \lambda_{n_m} =\mu$. 
Define 
\[ 
w_{nk}  = \left\{ \begin{array}{l}
 t,\,if\, k = l, \\
 1 - t,\,if\,k = m, \\
 0,\,  \textrm{otherwise,} \\
 \end{array} \right.
\]
 where $\ t \in \left( {0,1} \right)$. If there exist  $\delta>0$ and $K>0$ such that 
\[ 
\# \left\{ \lambda _{nj} ; {\left| {t\lambda _{nl}  + (1 - t)\lambda _{nm}  - \lambda _{nj} } \right|< \delta } \right\} < K,
\] 
then  $ \sigma_e(A) $ has a gap of width larger than $\ \delta$.\\
The advantage of this special choice is that we are able to avoid the assumptions on the bounds of $\sigma(A)$ and $\sigma_e(A)$. This shows that a  more general result is possible, provided that we choose the sequence of numbers $\ w_{nk} $ carefully.
    In the following theorem, we observe that the converse of Theorem (\ref{Kb1}) is true in the case of operators in the Arveson's class.
\bt Let A be a bounded self-adjoint operator in the Arveson's class. And suppose that there exists a gap in the essential spectrum. Then there exists a 
set of numbers $\ \{w_{nk} : k = 1,2,... n\}$ such that
 $\ 0 \leq w_{nk}\leq 1$ and $\ \sum\limits_{{\ k = 1}}^{\  n} {w_{nk}  = 1}$ and a $\ \delta >0$ such that
 \begin{equation} \nonumber
 \# \left\{ {\lambda _{nj} ;\left| {\sum\limits_{k = 1}^n {w_{nk} \lambda _{_{nk} }  - \lambda _{nj} } } \right| < \delta } \right\} < K,
 \end{equation}
 for some $ K>0$.\label{Kb2} \et
 \bp Let $\ \left( {a,b} \right)$ be a gap in the essential spectrum. Then by Theorem (\ref{Arv}), there exists sequences of eigenvalues of truncations $\ \lambda 
 _{n_l},\lambda _{n_m}$ such that
\begin{equation} \nonumber
\lim_{n_l\rightarrow\infty} \lambda_{n_l} = a \, \textrm{and} \lim_{n_m\rightarrow\infty} \lambda_{n_m} =b.
\end{equation}
Define
\begin{equation} \nonumber
 w_{nk}  = \left\{ \begin{array}{l}
 t,\,if\, k = l, \\
 1 - t,\,if\,k = m, \\
 0,\,\textrm{otherwise}, \\
 \end{array} \right.
 \end{equation}
 for some fixed $\ t \in \left( {0,1} \right)$. Since $c_t={ta+(1-t)b}\in \left( {a,b} \right)$, it is not an essential point. Also since $A$ is in the 
 Arveson's class, all such points are transient by Theorem (\ref{Arv1}).
Hence there exists a $\ \delta_1 >0$ such that sup$N_n \left( {c_t -\delta_1 ,c_t +\delta_1} \right)<K_1$ for some $K_1>0$. Also
\begin{equation} \nonumber
 \sum\limits_{{\  k = 1}}^{\ n} {w_{nk} \lambda _{nk} } = t\lambda_{n_l}+\left (1-t\right)\lambda_{n_m} {\rightarrow}{ta+(1-t)b}=c_t\,as\, n 
 {\rightarrow} \infty.
\end{equation}
 Therefore there exists an N such that
  $ \left| {c_t  - \sum\limits_{k = 1}^n {w_{nk} \lambda _{nk} } } \right| < {\raise0.7ex\hbox{$\delta_1 $} \!\mathord{\left/ {\vphantom {\delta_1  
  2}}\right.\kern-\nulldelimiterspace}
\!\lower0.7ex\hbox{$2$}} $ for all $n>N$.\\
Now if for some $n>N\,\,\,, \left| {\sum\limits_{k = 1}^n {w_{nk} \lambda _{_{nk} }  - \lambda _{nj} } } \right| < {\raise0.7ex\hbox{$\delta_1 $} 
\!\mathord{\left/ {\vphantom {\delta_1  2}}\right.\kern-\nulldelimiterspace}
\!\lower0.7ex\hbox{$2$}}$  then $ \left| {c_t  - \lambda _{nj} } \right| < \delta_1$.
Therefore,
\begin{equation} \nonumber
\# \left\{ {\lambda _{nj} ;\left| {\sum\limits_{k = 1}^n {w_{nk} \lambda _{_{nk} }  - \lambda _{nj} } } \right| < \frac{\delta_1}{2} } \right\} <N_n 
\left( {c_t -\delta_1 ,c_t +\delta_1} \right)<K_1, \, \forall \, n>N.
\end{equation}
 Now choosing $K=\textrm{sup}\{K_1,N\} \& \, \delta =\frac{\delta_1}{2} $, we complete the proof.
\ep
\br
In the above proof, numbers $\ \{w_{nk} : k = 1,2,... n\}$ and the bound $K$ will depend on the particular $t \in \left(0,1\right)$ that we choose.
\er
\subsection{Gaps under perturbation}
 Now we look at the spectral gaps that may occur between the bounds of the essential spectrum 
 of a holomorphic family of self-adjoint operators. Recall that 
the gaps remain invariant under a compact perturbation of the operator. 
The question that we address here is how stable these gaps, under a more general 
perturbation. Also the stability of the predictions of gaps under a holomorphic perturbation, 
is another question to be addressed here.

 We state the \textbf{stability theorem of bounded invertibility } and 
 use it to achieve some invariance for the gaps. The theorem is stated in a 
 more general form in \cite{Kat}. We need only the following special case.
 \bt \label{stability}
 Let A and B are bounded operators and A is invertible. If the quantity $\left\|A^{-1}\right\|\left\|B\right\|<1,$
  then A+B is also invertible.
 \et
The following theorem is an immediate consequence of the 
 stability theorem stated above. 
 \bt \label{perturbgap}
Let $\left(a,b\right)$ is a gap in $\sigma_e(A(0))$ which contains no discrete spectral value in it.
 Then for all small enough $\varepsilon>0$, 
 there exists a $\delta>0$ such that $\left(a+ \varepsilon ,b- \varepsilon \right)$ is a gap in the essential spectrum of the analytic family of operators
$A\left( x \right)\, \, \, \textrm{for every x with} \left|x\right|< \delta$.
\et
\bp
 First we note that, $A-\lambda I$ is 
 invertible for every $\lambda$ in the interval $\left(a,b\right),$ since it contains no spectral value.
 Therefore,
  \begin{equation} \nonumber
\textrm{sup}\left\{ {\left\| {\left( {A - \lambda I} \right)^{ - 1} } \right\|;\lambda  \in \left( {a + \varepsilon ,b - \varepsilon } \right)} 
\right\}=M < \infty \, \, \textrm{for a fixed}\,\, \varepsilon>0.
\end{equation}
Now using the continuity assumption, corresponding to minimum of $\{\frac{1}{M},\epsilon\},$ there exists a $\delta>0,$ such that
\begin{equation} \nonumber
\left\|(A(x)-A(0)\right\|< \min \{\frac{1}{M},\epsilon\}\, \, \, \textrm{for every x with} \left|x\right|< \delta.
\end{equation}
Now for $\left|x\right|< \delta,$ observe that
\begin{equation} \nonumber
\left\|(A-\lambda I)^{-1}\right\|\left\|A(x)-A(0)\right\|<M.\frac{1}{M}<1
\end{equation}
for every $\lambda$ in the interval 
$\left(a+ \varepsilon ,b- \varepsilon \right)$. 

Hence by Theorem (\ref{stability}), if $\left|x\right|< \delta,$ then 
\[A(x)-\lambda I=A(x)-A(0)+A(0)-\lambda I\] is invertible for every $\lambda$ in the interval 
$\left(a+ \varepsilon ,b- \varepsilon \right)$.
Therefore the interval $\left(a+ \varepsilon ,b- \varepsilon \right)$ does not intersect with
$\sigma(A\left( x \right)),\,\,\textrm{for every x with} \left|x\right|< \delta$. 

Now, since $\left\|A(x)-A(0)\right\|<\epsilon,$ 
$\left(a+ \varepsilon ,b- \varepsilon \right)$ will lie between the bounds of $\sigma_e(A(x))$, 
for every x with $\left|x\right|< \delta$.
We conclude that $\left(a+ \varepsilon ,b- \varepsilon \right)$ is a spectral gap in 
$\sigma_e(A(x))$ for all x, with $\left|x\right|<\delta$.
\ep
\br
In Theorem (\ref{perturbgap}), $\epsilon$ must be small enough so that the interval 
$\left(a+ \varepsilon ,b- \varepsilon \right)$ should makes sense. This theorem indicates that to some extend,
the gaps are stable under small norm perturbation. Once we get $\left(a+ \varepsilon ,b- \varepsilon \right)$ 
is a gap, we may remove that interval and look at the rest of the interval $\left(a,b\right)$ and continue the 
search for gaps.
\er
Let's look at an example to support the above theorem.
\be \label{borgpert}
Define a two parameter family of matrix valued symbols as follows

\begin{equation}\label{fx} \nonumber
f(x, \theta)= \left[ \begin{array}{cccccccccc}
    a_1(x) & 1 & & & & e^{-i\theta}\\
      1 & a_2(x)  & 1\\
      & 1 & a_3(x)& 1\\
       &&1 & a_4(x) & 1\\
      & && \ddots &\ddots& \ddots  \\
      e^{i\theta}& & &    &1 & a_p(x)\\
    \end{array} \right],
\end{equation}
where $a_1(.),a_2(.) \ldots a_p(.)$ are real analytic functions and $\theta$ varying in the interval $[0,2\pi]$. 
Note that 
\begin{equation}\nonumber
f(x, \theta)=A_0(x)+A_1 e^{i\theta}+A_{-1} e^{-i\theta},\,\,\textrm{where}
\end{equation}
\[
 A_0(x)=\left[ {\begin{array}{*{20}c}
   {a_1(x) } & 1 & {} & {} & {} & { }  \\
   1 & {a_2(x)} & 1 & {} & {} & {}  \\
   {} & 1 & . & . & {} & {}  \\
   {} & {} & . & . & . & {}  \\
   {} & {} & {} & . & . & 1  \\
   {}  & {} & {} & {} & 1 & {a_p(x) }  \\
\end{array}} \right],
\]
\[
 A_1=\left[ {\begin{array}{*{20}c}
   { } &  & {} & {} & {} & { }  \\
    & {} &  & {} & {} & {}  \\
   {} &  &  &  & {} & {}  \\
   {} & {} &  &  &  & {}  \\
   {} & {} & {} &  &  &   \\
   {1}  & {} & {} & {} &  & {}  \\
\end{array}} \right]= {A_{-1}}^T.
\] 
We consider the one parameter family of block Toeplitz-Laurent operators arising from these symbols, which are represented 
by the following doubly infinite matrices,
\[
 A(x) = \left[ {\begin{array}{*{20}c}
   {\ddots } & {\ddots } & {} & {} & {} & {} & {} & {} & {} & {}  \\
   {\ddots} & {A_0(x) } & {A_{ - 1} } & {} & {} & {} & {} & {} & {} & {}  \\
   {} & {A_1 } & {A_0(x) } & {A_{ - 1} } & {} & {} & {} & {} & {} & {}  \\
   {} & {} & A_1 & A_0(x) &A_{ - 1} & {} & {} & {} & {} & {}  \\
   {} & {} & {} & A_1 & A_0(x) & A_{ - 1} & {} & {} & {} & {}  \\
   {} & {} & {} & {} & A_1 & A_0(x) &A_{ - 1} & {} & {}   \\
  {} & {} & {} & {} & {} & {A_1 } & {A_0(x) } & {\ddots }  \\
{} & {} & {} & {} & {} & {} & {\ddots } & {\ddots}  \\
\end{array}} \right]
\]

 Thus we get an 
analytic family of bounded operators, $A(x)$ which are self-adjoint for all $x$ in the domain.
Now, by Borg's theorem for discrete Schrodinger operator \cite{Fla,KSN}, the essential spectrum of 
$A(x_0)$, has no gaps if and only if $a_1(x_0)=a_2(x_0) \ldots =a_p(x_0)$. Hence if there is a gap in 
$\sigma_e(A(0))$ then $a_i(0)<a_{i+1}(0)$ for 
some $i$. Using the continuity of $a_i\, \&\, a_{i+1}$, we can find a $\delta>0$ such that  $a_i(x)<a_{i+1}(x)$
for all x with $\left|x\right|<\delta$. Hence there is a gap for $A(x)$ for all such $x$.
\ee

Using Theorem (\ref{Kb1}) and Theorem (\ref{perturbgap}), we arrive at the following conclusions.
The gap predictions that we have done for a single operator, are remain valid for a family of operators. The advantage is that we can predict gaps of a family of operators, with 
assumptions only on the unperturbed operator. We give the precise statement below.
\bc
 Let A(x) be a holomorphic family of operators with
 $A(0)=A$, and $ \lambda_{n1}(A_n)\geq\lambda_{n2}(A_n)\geq...\geq\lambda_{nn}(A_n)$ be the eigenvalues 
of $A_n$ arranged in decreasing order. For each  positive integer n, let $\ \{w_{nk} : k = 1,2,... n\}$  be a set of numbers such that $\ 0 \leq 
w_{nk}\leq 1$ and
 $ \sum\limits_{k = 1}^n w_{nk}=1$ and suppose there exists a $\delta>0$ and $ K>0$ such that
\begin{equation}\nonumber
\# \left\{ {\lambda _{nj} ;\left| {\sum\limits_{k = 1}^n {w_{nk} \lambda _{_{nk} }  - \lambda _{nj} } } \right| < \delta } \right\} < K.
\end{equation}
In addition, if we assume that $\ \sigma_e(A) $ and $\ \sigma(A)$ coincide, then $\ \sigma_e(A(x))$
 has a gap for each x in a sufficiently small neighborhood of 0.
\ec
\bp
By Theorem (\ref{Kb1}), A(0) has a spectral gap. By Theorem (\ref{perturbgap}), 
there exists a neighborhood of 0, with $A(x)$ has gaps for all x in the neighborhood. Hence the proof. 
\ep
\br
In the above case, we considered perturbation of operators and not the perturbations of their truncations. In the 
Example (\ref{borgpert}) also the perturbed symbol is directly related to the perturbation of operators. The perturbation of truncations and their link with the spectrum of the original operator is another problem yet to be handled.  
\er
\section{Gap issues of Block Toeplitz-Laurent operators}
In this section, we look at the spectral gap issues of some block Toeplitz-Laurent operators. The operators under our concern are some perturbations of discrete Schrodinger operator on $l^2(Z)$. Below we try to improve the discrete version of Borg's theorem in a more general set up with the linear algebraic techniques used in \cite{KSN}. We refer to \cite{KSN} for a detailed description of the discretization of Schrodinger operator and formulating the discrete version in terms of block Toeplitz-Laurent operators and the matrix valued symbol.
\bt \label{Borg}
Let A be the bounded operator defined by the block Toeplitz-Laurent matrix
\[
 A = \left[ {\begin{array}{*{20}c}
   {\ddots } & {\ddots } & {} & {} & {} & {} & {} & {} & {} & {}  \\
   {\ddots} & {A_0 } & {A_{ - 1} } & {A_{ - 2}} & {} & {\ldots} & {A_{ - N}} & {\ldots} & {} & {}  \\
   {} & {A_1 } & {A_0 } & {A_{ - 1} } & {A_{ - 2}} & {} & {\ldots} & {A_{ - N}} & {\ldots} & {}  \\
   {} & {A_{ 2}} & A_1 & A_0 &A_{ - 1} & {A_{ - 2}} & {} & {\ldots} & {A_{ - N}} & {\ldots}  \\
   {} & {} & {A_{ 2}} & A_1 & A_0 & A_{ - 1} & {A_{ - 2}} & {} & {\ldots} & {A_{ - N}}  \\
   {\ldots} & {A_{ N}} & {\ldots} & {A_{ 2}} & A_1 & A_0 &A_{ - 1} & {A_{ - 2}} & {}   \\
   {} & {\ldots} & {A_{ N}} & {\ldots} & {A_{ 2}} & A_1 & A_0 & A_{ - 1} & {A_{ - 2}} & {}  \\
   {} & {} & {\ldots} & {A_{ N}} & {\ldots} & {A_{ 2}} & A_1 & A_0 & A_{ - 1} & {A_{ - 2}}  \\
   {} & {} & {} & {\ldots} & {A_{ N}} & {\ldots} & {A_{ 2}} & {A_1 } & {A_0 } & {\ddots }  \\
   {} & {} & {} & {} & {} & {} & {} & {} & {\ddots } & {\ddots}  \\
\end{array}} \right]
\]
where 
\[
 A_0  = \left[ {\begin{array}{*{20}c}
   {b_1 } & 1 & {} & {} & {} & {a_0}  \\
   1 & {b_2 } & 1 & {} & {} & {}  \\
   {} & 1 & . & . & {} & {}  \\
   {} & {} & . & . & . & {}  \\
   {} & {} & {} & . & . & 1  \\
   {a_0} & {} & {} & {} & 1 & {b_p }  \\
\end{array}} \right], \ \ \
 A_{ k}  = \left[ {\begin{array}{*{20}c}
   {} & {} & {} & {} & {} & {a_k}  \\
   {} & {} & {} & {} & {} & {}  \\
   {} & {} & {} & {} & {} & {}  \\
   {} & {} & {} & {} & {} & {}  \\
   {} & {} & {} & {} & {} & {}  \\
   {} & {} & {} & {} & {} & {}  \\
\end{array}} \right]={A_{ -k}}^T,
\]
such that $b_1\leq b_2 \ldots \leq b_p$ and $\sum_k \left|a_k\right|< \infty$.
If A has connected essential spectrum then $b_1=b_2 \ldots =b_p$.
\et
\bp
The matrix-valued symbol associated with the block Toeplitz-Laurent operator A is
 \[
 \tilde{f}\left( \theta  \right) = \left[ {\begin{array}{*{20}c}
   {b_1 } & 1 & {} & {} & {} & {f\left(\theta \right) }  \\
   1 & {b_2 } & 1 & {} & {} & {}  \\
   {} & 1 & . & . & {} & {}  \\
   {} & {} & . & . & . & {}  \\
   {} & {} & {} & . & . & 1  \\
   {\bar{f\left(\theta \right)} } & {} & {} & {} & 1 & {b_p }  \\
\end{array}} \right].
\]
where $f\left(\theta \right)= \sum_k a_ke^{ ik\theta }$. Therefore from \cite{bottcher} we have
\begin{equation} \label{bottcher-id}
\sigma_{\hbox{ess}}(A)= \bigcup_{j=1}^p \left[\inf_\theta(\lambda_j(\tilde{f}\left( \theta  \right)),\sup_\theta(\lambda_j(\tilde{f}\left( \theta  \right))\right].
\end{equation}
Now consider the sub matrices 
\[
\begin{array}{l}
 P_1  = \left[ {\begin{array}{*{20}c}
   {b_1 } & 1 & {} & {} & {} & {}  \\
   1 & {b_2 } & 1 & {} & {} & {}  \\
   {} & 1 & . & . & {} & {}  \\
   {} & {} & . & . & . & {}  \\
   {} & {} & {} & . & . & 1  \\
   {} & {} & {} & {} & 1 & {b_{p-1} }  \\
\end{array}} \right],  
 \,\,\,\,\,\,\,\,\,\,\, 
 P_2  = \left[ {\begin{array}{*{20}c}
   {b_2 } & 1 & {} & {} & {} & {}  \\
   1 & {b_3 } & 1 & {} & {} & {}  \\
   {} & 1 & . & . & {} & {}  \\
   {} & {} & . & . & . & {}  \\
   {} & {} & {} & . & . & 1  \\
   {} & {} & {} & {} & 1 & {b_p }  \\
\end{array}} \right]. \\ 
 \end{array}
\]
If any of their eigenvalues are different, say $\lambda_j(P_1)<\lambda_j(P_2)$, then by Cauchy Interlacing theorem,
$\lambda_{j}(\tilde{f}\left( \theta  \right))\leq\lambda_j(P_1)<\lambda_j(P_2)\leq \lambda_{j+1}(\tilde{f}\left( \theta  \right)), \,\,\forall\,\, \theta$.
But from (\ref{bottcher-id}), this will give us the contradiction that essential spectrum of A is not connected.
 Hence all the eigenvalues of $P_1$ and $P_2$ are same. Therefore
\begin{equation}\nonumber
\hbox{trace}(P_1)-\hbox{trace}(P_2)=b_1-b_p=0
\end{equation}
Hence $b_1=b_2 \ldots =b_p$.
\ep
\br
The converse of the above 
assertion is in general not true. We may have gaps even if the diagonal entries of the
block Toeplitz-Laurent operator are same.
For if A is the block Toeplitz-Laurent operator arising from the matrix valued symbol
 \[
 \tilde{f}\left( \theta  \right) = \left[ {\begin{array}{*{20}c}
   {b} & {1+f\left(\theta \right) }  \\
   {1+\bar{f\left(\theta \right)} } & {b } \\
\end{array}} \right].
\]
where f is a non negative function, then the eigenvalue functions of $\tilde{f}\left( \theta  \right)$ are 
\begin{equation} \nonumber
\lambda_1(\theta)= b-1-f(\theta), \,
\lambda_2(\theta)= b+1+f(\theta)
\end{equation}
Hence spectrum of A will have a gap, since f is non negative.
\er
\br
We remark that the diagonal entries correspond to the periodic potential of the discrete Schrodinger operator. Hence we have proved the discrete Borg-type theorem for a perturbed operator with some extra assumptions on the potential.
\er
\be
The assumption $b_1\leq b_2 \ldots \leq b_p$ can not be dropped in the above theorem, if $p>2$. For if we
consider the block Toeplitz-Laurent operator arising from the matrix valued symbol 
\[
 \tilde{f}\left( \theta  \right) = \left[ {\begin{array}{*{20}c}
   {1} & {1}&{0}&{10cos\left(\theta \right) }  \\
   {1} & {2}&{1}&  {0} \\
   {0} & {1}&{2}&{1} \\
   {10cos\left(\theta \right)}  &{0}& {1}&{1 } \\
\end{array}} \right].
\]
The eigenvalue functions of $\tilde{f}\left( \theta  \right)$ are 
\begin{equation} \nonumber
\lambda_{1,2}(\theta)= 2+5cos(\theta) \pm \sqrt{25cos^2(\theta)-10cos(\theta)+2}
\end{equation}
\begin{equation} \nonumber
\lambda_{3,4}(\theta)= 1-5cos(\theta) \pm \sqrt{25cos^2(\theta)+1}
\end{equation}
We list the values of these functions at certain points in the table below.
\begin{table}[h]%
\begin{tabular}{|c|c|c|c|c|}
\hline
$\theta$ & $\lambda_1(\theta)$&$\lambda_2(\theta)$&$\lambda_3(\theta)$&$\lambda_4(\theta)$\\
\hline
0&11.123&2.877&1.099&-9.099\\
\hline
$\pi$&3.083&-9.083&11.099&.901\\
\hline
\end{tabular}
\end{table}

From the table, it is clear that the ranges of the above continuous functions intersect. Hence their union is a 
connected interval. Therefore the essential spectrum of the operator has no gaps, even the periodic potential does not
reduce to a constant.
\ee
\subsection{Perturbation of matrices}
Finally we use some known results on the 
 bounds for the eigenvalues of perturbed matrices (see \cite{Bha1},\cite{Li} and references there in)
 to strengthen our results by viewing the matrix valued symbol as a 
perturbation of some constant matrix. 
\bl \label{matrixpert}
Let $ H = \left( {\begin{array}{*{20}c}
   {H_1 } & E  \\
   {E^* } & {H_1 }  \\
\end{array}} \right) $ and
$ \tilde{H} = \left( {\begin{array}{*{20}c}
   {H_1 } & 0 \\
   {0 } & {H_1 }  \\
\end{array}} \right) $, $\lambda_1 \geq\lambda_2 \geq \ldots \lambda_p$ and $\tilde{\lambda_1} \geq\tilde{\lambda_2} \geq \ldots \tilde{\lambda_p}$
be the eigenvalues respectively. Then
\begin{equation}\label{eigenvalue estimate 1}
\left| {\lambda _j  - \tilde{\lambda_j}} \right| \le \left\| E \right\|
\end{equation}
\el
\bt \label{perturbborg}
Let A be the operator considered in Theorem (\ref{Borg}). If $\lambda_1 \geq\lambda_2 \geq \ldots \lambda_p$ are
eigenvalues of the matrix 
\[
 \left[ {\begin{array}{*{20}c}
   {b_1 } & 1 & {} & {} & {} & { }  \\
   1 & {b_2 } & 1 & {} & {} & {}  \\
   {} & 1 & . & . & {} & {}  \\
   {} & {} & . & . & . & {}  \\
   {} & {} & {} & . & . & 1  \\
   {}  & {} & {} & {} & 1 & {b_p }  \\
\end{array}} \right],
\]
and $\sigma_{\hbox{ess}}(A)$ has no gap, then $\left| {\lambda _j  - \lambda _{j+1} } \right| \le 2\left\| f \right\|_\infty \, \forall j=1,2 \ldots p-1.$
In addition, if we assume that $b_1\leq b_2 \ldots \leq b_p$, then 
$\left| {\lambda _j  - \lambda _j } \right| \le 2\, \forall j=1,2 \ldots p.$
\et
\bp
Apply above lemma with 
 \[
 H(.)= \left[ {\begin{array}{*{20}c}
   {b_1 } & 1 & {} & {} & {} & {f\left(\theta \right) }  \\
   1 & {b_2 } & 1 & {} & {} & {}  \\
   {} & 1 & . & . & {} & {}  \\
   {} & {} & . & . & . & {}  \\
   {} & {} & {} & . & . & 1  \\
   {\bar{f\left(\theta \right)} } & {} & {} & {} & 1 & {b_p }  \\
\end{array}} \right],
\tilde{H}=\left[ {\begin{array}{*{20}c}
   {b_1 } & 1 & {} & {} & {} & { }  \\
   1 & {b_2 } & 1 & {} & {} & {}  \\
   {} & 1 & . & . & {} & {}  \\
   {} & {} & . & . & . & {}  \\
   {} & {} & {} & . & . & 1  \\
   {}  & {} & {} & {} & 1 & {b_p }  \\
\end{array}} \right] 
\]
and 
\[
E= \left[ {\begin{array}{*{20}c}
   { } &  & {} & {} & {} & {f\left(\theta \right) }  \\
   {} & {} & {} & {} & {} & {}  \\
   {} & {} &  & {} & {} & {}  \\
   {} & {} & {} & {} & {} & {}  \\
   {} & {} & {} & {} & {}& {}  \\
   {}  & {} & {} & {} & {} & { }  \\
\end{array}} \right].
\]
Then we get 
\begin{equation} \nonumber
\left| {\lambda _j  - \tilde{\lambda_j(\theta)}} \right| \le \left\| E \right\|=\left\| f \right\|_\infty
\end{equation} by (\ref{eigenvalue estimate 1}). Combining with (\ref{bottcher-id}), we get
\begin{equation}\nonumber
\sigma_{\hbox{ess}}(A)= \bigcup_{j=1}^p \left[\inf_\theta(\lambda_j(f_s(\theta))),\sup_\theta(\lambda_j(f_s(\theta)))\right] \subseteq 
\bigcup_{j=1}^p \left[\lambda_j-\left\| f \right\|_\infty,\lambda_j+\left\| f \right\|_\infty\right].
\end{equation}
Therefore if $\left| {\lambda _j  - \lambda _{j+1} } \right| > 2\left\| f \right\|_\infty\, \textrm{for some}\,\, j$, then there exists a gap in the essential spectrum. Hence we proved the first assertion. Now in addition, if we assume that $b_1\leq b_2 \ldots \leq b_p$, then 
since the essential spectrum of A is connected, by theorem (\ref{Borg}), $b_1=b_2 \ldots =b_p$. This implies that H is a tridiagonal Toeplitz matrix with $b_1$ on diagonal and 1 as off diagonal entry. The eigenvalues of such matrices are explicitly known and they are $b_1+2cos(\frac{\pi k}{p+1})$. So the second assertion follows by a simple computation.
\ep
\br
Using the last theorem, we can predict the nature of spectrum of the operator A with only looking at the eigenvalues of the $p\times p$ matrix $\tilde{H}.$
\er
\subsection{Jacobi matrices}
We can apply the same technique in the case of periodic Jacobi matrices to predict gaps in the
essential spectrum.

Define a double infinite, p-periodic, $p\geq2,$ real Jacobi matrix by
\begin{equation}\label{jacobi1}
J= \left[ \begin{array}{cccccccccc}
 \ddots & \ddots & \\
\ddots  & b_1 &  a_1 & \\
  & a_1&b_2& a_2 \\
  &    &\ddots & \ddots & \ddots\\
 &     &      & \ddots & b_{p} & a_p \\
   &     &      &    & a_p & b_1 & \ddots \\
 &    &     &    &  & \ddots 	& \ddots
     \end{array} \right],a_{n+p}=a_n > 0; b_{n+p} = b_n:
\end{equation}
 following the standard convention $a_n >0.$ 
  An important observation is that J is the block Toeplitz-Laurent operator, where in the case $p\geq3,$ the
symbols are given by
\[
f_k\left( \theta  \right) = \left[ {\begin{array}{*{20}c}
   {b_{k+1} } &a_{k+1} & {0} & {} & {} & {e^{i\theta }a_{k+p-1} }  \\
   a_{k+1} & {b_{k+2}} & a_{k+2} & {} & {} & {0}  \\
   {0} & a_{k+2} & . & . & {} & {}  \\
   {} & {} & . & . & . & {}  \\
   {} & {} & {} & . & . & a_{k+p-1}  \\
   {e^{ - i\theta }a_{k+p-1} } & {0} & {} & {} & a_{k+p-1} & {b_{k+p} }  \\
\end{array}} \right],k=0,1,\ldots p-1,
\]
Also the spectrum of J is given by the following identity.
\begin{equation}\label{Jacobispectrum}
\sigma(J)=\bigcup_{j=1}^p \left[\min_{\theta}\lambda_j(f_k(\theta)),\max_{\theta}\lambda_j(f_k(\theta))\right].
\end{equation}
\bc Let J be the Jacobi matrix defined by (\ref{jacobi1}). If $\lambda_1 \geq\lambda_2 \geq \ldots \lambda_p$ are
eigenvalues of the matrix 
\[
 \left[ {\begin{array}{*{20}c}
   {b_1 } & {a_1} & {} & {} & {} & { }  \\
   {a_1} & {b_2 } & {a_2} & {} & {} & {}  \\
   {} & {a_2} & . & . & {} & {}  \\
   {} & {} & . & . & . & {}  \\
   {} & {} & {} & . & . & {a_{p-1}}  \\
   {}  & {} & {} & {} & {a_{p-1}} & {b_p }  \\
\end{array}} \right].
\]
Then $\sigma_{\hbox{ess}}(J)$ has a gap, if $\left| {\lambda _j  - \lambda _{j+1} } \right| > 2\left| a_{p-1} \right|$ 
for some j.
\ec
\bp
The proof is an imitation of the proof of Theorem (\ref{perturbborg}), however all the 
details are provided here. Apply Lemma (\ref{matrixpert}) with 
 \[
 H(.)= \left[ {\begin{array}{*{20}c}
   {b_1 } & {a_1} & {} & {} & {} & {a_{p-1}e^{i\theta}}  \\
   {a_1} & {b_2 } & {a_2} & {} & {} & {}  \\
   {} & {a_2} & . & . & {} & {}  \\
   {} & {} & . & . & . & {}  \\
   {} & {} & {} & . & . & {a_{p-1}}  \\
   {{a_{p-1}e^{-i\theta} }}  & {} & {} & {} & {a_{p-1}} & {b_p }  \\
\end{array}} \right],
\tilde{H}=\left[ {\begin{array}{*{20}c}
   {b_1 } & {a_1} & {} & {} & {} & { }  \\
   {a_1} & {b_2 } & {a_2} & {} & {} & {}  \\
   {} & {a_2} & . & . & {} & {}  \\
   {} & {} & . & . & . & {}  \\
   {} & {} & {} & . & . & {a_{p-1}}  \\
   {}  & {} & {} & {} & {a_{p-1}} & {b_p }  \\
\end{array}} \right],
\]
and 
\[
E= \left[ {\begin{array}{*{20}c}
   { } &  & {} & {} & {} & {a_{p-1}e^{i\theta}}  \\
   {} & {} & {} & {} & {} & {}  \\
   {} & {} &  & {} & {} & {}  \\
   {} & {} & {} & {} & {} & {}  \\
   {} & {} & {} & {} & {}& {}  \\
   {}  & {} & {} & {} & {} & { }  \\
\end{array}} \right].
\]
Then we get 
\begin{equation} \nonumber
\left| {\lambda _j  -\lambda_j(H(\theta)} \right| \le \left\| E \right\|=\left\|{a_{p-1}e^{i\theta}}\right\|_\infty=
\left|{a_{p-1}}\right|\,\,
\textrm{by (\ref{eigenvalue estimate 1})}.
\end{equation} Combining with (\ref{Jacobispectrum}), we get
\begin{equation}\nonumber
\sigma_{\hbox{ess}}(J)= \bigcup_{j=1}^p \left[\inf_\theta(\lambda_j(H(\theta))),\sup_\theta(\lambda_j(H(\theta)))\right] \subseteq 
\bigcup_{j=1}^p \left[\lambda_j-\left|{a_{p-1}}\right|,\lambda_j+\left|{a_{p-1}}\right|\right].
\end{equation}
Therefore if $\left| {\lambda _j  - \lambda _{j+1} } \right| > 2\left|{a_{p-1}}\right|\, \textrm{for some}\,\, j$, then there exists a gap in the essential spectrum. Hence the proof.
\ep
\br
The last couple of theorems help us to reduce the computations in predicting spectral gaps, 
for operators arising from the matrix valued symbols.
We need to check only the eigenvalues of a matrix with constant entries. 
The proof also gives us the spectral inclusion 
\begin{equation}\nonumber
\sigma_{\hbox{ess}}(A) \subseteq 
\bigcup_{j=1}^p \left[\lambda_j-\left\| f \right\|_\infty,\lambda_j+\left\| f \right\|_\infty\right].
\end{equation}
which is very important, since the right hand side includes only the eigenvalues of a constant matrix.
Whether equality holds in this inclusion, is still not clear to us.
\er
\section{Concluding Remarks}
  We conclude this note by listing down some remarks and future problems.
\begin{itemize}
\item Using Theorem (\ref{Kb1}) and the special choice I, we could predict the existence of spectral gaps from the finite matrix entries. Theorem (\ref{perturbborg}) and its corollary can be used to predict the spectral gaps of the corresponding operators, by looking at the eigenvalues of a finite matrix with constant entries.
\item  The Borg's theorem is a classical theorem in inverse spectral theory. The discrete versions are also folklore (\cite{Fla}). The techniques of the proof here, are adapted from \cite{KSN}.
\item     
 The discrete spectral values lying between a gap in the essential spectrum, can be computed using linear algebraic techniques.
To see this, let $ \left( {a,b} \right)$ be a gap in the essential spectrum of A. Let $ \lambda _0  = (a+b)/2 $. Since$\ \lambda _0$ is in the gap, $ 
f\left( {\lambda _0 } \right)>0$, where $ f\left( {\lambda _0 } \right)$ is the lower bound of the essential spectrum of
$ {\left( {A - \lambda _0 I} \right)^2 }$, all the discrete spectral values below that can be computed with the use of truncations by Theorem (\ref{Bot}). If
$\beta$ is an eigenvalue in the gap, $ \left( {\beta  - \lambda_0 } \right)^2 $ will be an eigenvalue lying below the lower bound of the essential 
spectrum of $\left( {A - \lambda _0 I} \right)^2 $. From these we can compute $ \beta$.

Looking at these observations under a holomorphic perturbation is an interesting problem.
       \item Also under compact perturbation, though the spectral gaps remain the same, discrete eigenvalues may appear or disappear inside such gaps. 
       Another problem is to handle such situations linear algebraically.
       \item Another scope is to carry over these results to the case of unbounded operators. In particular one may think of estimating the spectrum 
       and spectral gaps of Schrodinger operators by the eigenvalues of its truncations.
\end{itemize}
\textbf{Acknowledgments:} Kiran Kumar and M.N.N. Namboodiri would like to thank Prof. Yuri Safarov, Department of Mathematics,
King's College London, for useful discussions and suggestions during his visit to Kerala under the Erudite program of Govt. of Kerala.
Kiran Kumar is thankful to CSIR, KSCSTE for financial support. Stefano Serra-Capizzano is thankful to the Govt. of Kerala, Erudite program, to the 
Italian MiUR, PRIN 2008 N. 20083KLJEZ, for financial support.

\end{document}